%% file: 2020-Lyap.tex
\documentclass[letterpaper,10pt, conference,twoside]{ieeeconf}  

\usepackage{amsmath,amssymb,color,cite}
\usepackage{graphicx}

\usepackage{latexsym}
\usepackage{verbatim}
\usepackage{cite}
\usepackage[english]{babel}
\usepackage[latin1]{inputenc}
\usepackage{enumerate}

\usepackage{amsthm}

\usepackage{tikz}
\usetikzlibrary{tikzmark}
\usepackage{amsfonts}
\usepackage{amsbsy}
\usepackage{bm}

\let\leq\leqslant
\let\geq\geqslant
\let\emptyset\varnothing

\input{ka-newcommands-siam}

\newcounter{todocounter}
\usepackage[colorinlistoftodos]{todonotes}

\newtheorem{theorem}{Theorem}
\newtheorem{lemma}[theorem]{Lemma}
\newtheorem{cor}[theorem]{Corollary}
\theoremstyle{definition}
\newtheorem{definition}[theorem]{Definition}
\newtheorem{example}[theorem]{Example}

\newtheorem{conjecture}[theorem]{Conjecture}
\newtheorem{remark}[theorem]{Remark}

\newcommand{\ico}{\scalebox{0.6}{ $\square$ }}

\title{On Duality for Lyapunov Functions of Nonstrict Convex Processes}

\author{Jaap Eising, M. Kanat Camlibel}

\begin{document}

\maketitle

\renewcommand{\thefootnote}{\fnsymbol{footnote}}

\footnotetext{The authors are with the Bernoulli Institute for Mathematics, Computer Science, and Artificial Intelligence, University of Groningen, Nijenborgh 9, 9747 AG, Groningen, The Netherlands. (email: {\tt j.eising@rug.nl; m.k.camlibel@rug.nl})}

\begin{abstract}
This paper provides a novel definition for Lyapunov functions for difference inclusions defined by convex processes. It is shown that this definition reflects stability properties of nonstrict convex processes better than previously used definitions. In addition the paper presents conditions under which a weak Lyapunov function for a convex process yields a strong Lyapunov function for the dual of the convex process. 
\end{abstract}

\section{Introduction}
In this paper we study Lyapunov functions for systems given by difference inclusions of the form 
\[ x_{k+1} \in H(x_k),\]
where $H:\mathbb{R}^n\rightrightarrows \mathbb{R}^n$ is a convex process: a set-valued map whose graph is a convex cone. Since being introduced in \cite{Rockafellar:67,Rockafellar:70} these maps have attracted attention for a few reasons. Mainly, these systems can used to describe the dynamics of conically constrained linear systems. Such conic constraints encompass for example nonnegative controls, systems evolving only in the nonnegative orthant, and more general dynamics arising from physical modeling. This means that systems of this form are encountered in various applications, including von Neumann-Gale economic growth models \cite{makarov:77}, cable-suspended robots \cite{j37,oh:05-2} and chemical reaction networks \cite{Angeli:09}. Lastly, as shown in e.g. \cite{Frankowska:87,Aubin:84}, difference inclusions of convex processes can be used as approximations of more complex set-valued maps. As such, local properties of more general difference inclusions can be described in terms of properties of an approximating convex process. 

All of these reasons have led to interest in the analysis of systems described by convex processes. In \cite{AFO:86} and \cite{PD:94} the controllability problem is resolved for strict (nonempty everywhere) convex processes in continuous-time and discrete-time, respectively. The works \cite{Smirnov:02} and \cite{Phat:96} similarly characterize stabilizability for strict convex processes. An important ingredient of these results is the use of duality: Controllability and stabilizability of a convex process can be characterized in terms of spectral properties of the dual (or adjoint) convex process. 

However, the assumption of strictness excludes many interesting dynamical systems from consideration. In terms of constrained linear systems, this assumption allows for only constraints on the inputs, but not on the states. Moreover, as shown in \cite{Seeger:01}, the aforementioned duality relation cannot be established for nonstrict convex processes in general. This has led to the work in \cite{j37} and \cite{ReachNullc:19}: These papers give conditions on the domain of a convex process under which the previously mentioned results can be generalized.

A similar duality relation arises when considering Lyapunov functions for difference inclusions of convex processes. As is the case for linear systems, Lyapunov functions have proven to be a natural and useful tool in stability analysis of difference inclusions, as evidenced by e.g. \cite{Kellet:04a ,Kellet:04 ,Angeli:09 ,Goebel:11}. For the more specific class of convex processes, \cite{Smirnov:02} already employed Lyapunov functions in proving its relation between stabilizability of a primal system and stability of the dual. This relation was made much more explicit in the extensive treatise on Lyapunov functions given in \cite{Goebel:13}. A result of particular interest is Theorem 2.2 of \cite{Goebel:13}. This theorem reveals for strict convex processes that the convex conjugate of a weak Lyapunov function for the primal system is a (strong) Lyapunov function for the dual. 

This last result is the starting point of the investigation in this paper. We study weaker conditions on the domain of a convex process under which a similar result can be shown to hold. Our contribution is twofold: First, we modify the definitions of weak and strong Lyapunov functions to better reflect the stability properties of nonstrict convex processes. We will illustrate this by means of a few examples. Second, we will prove results generalizing those of, among others, Theorem 2.2 of \cite{Goebel:13}. To be precise, we show that we can obtain a strong Lyapunov function for the dual by taking a restriction of a weak Lyapunov function of the primal and then taking the convex conjugate. The domain conditions required to obtain this result  can be shown to be in some sense close to those required for the existence of such a strong Lyapunov function. 

This paper is organized as follows: We start with some preliminary knowledge on convex processes and extended real-valued functions in Section~\ref{sec:cp} and Section~\ref{sec:ervf} respectively. After this, we define weak and strong Lyapunov functions and motivate this definition in Section~\ref{sec:Lyap}. In Section~\ref{sec:results} we provide the main results of the paper and lastly we provide conclusions in Section~\ref{sec:conc}.

\section{Convex processes}\label{sec:cp}
Given a convex set $\calS\subseteq\mathbb{R}^n$, we denote its closure by $\cl\calS$ and its relative interior by $\rint\calS$. Given in addition the convex set $\calT$ and scalar $\rho\in\R$ we define the sum and scalar product of sets as: 
\[ \calS+\calT = \{s+t \mid s\in\calS, t\in\calT\}, \quad \rho \calS = \{\rho s \mid s\in\calS\}. \] 
We consider \textit{convex cones}: nonempty convex sets that are closed under nonnegative scalar multiplication. 

A set-valued map $H:\R^n\rightrightarrows \R^n$ is called a {\em convex process\/}, a {\em linear process\/}, {\em closed\/} if its graph
\[ \graph(H) = \{(x,y)\in \R^n\times \R^n \mid y\in H(x)\}\]
is a convex cone, a subspace, closed, respectively.

The \textit{domain} and \textit{image} of $H$ are defined as $\dom(H)= \{x\in \R^n\mid H(x)\neq \emptyset \}$ and $\im (H) = \{y\in \R^n : \exists\, x \textrm{ s.t. } y\in H(x)\}$. If $\dom(H)=\R^n$, we say that $H$ is \textit{strict}.

In this paper we consider systems described by a \textit{difference inclusion} of the form:
\begin{equation}\label{eq:diffinc} x_{k+1} \in H(x_k) \quad k\geq 0,
\end{equation}
where $H:\R^n\rightrightarrows \R^n$ is a convex process. By a \textit{trajectory} of \eqref{eq:diffinc}, we mean a sequence $(x_k)_{k\in\N }$ such that \eqref{eq:diffinc} holds for all $k\geq 0$. The \textit{behavior} (see e.g. \cite{Willems:91}) is the set of all trajectories: 
\[ \mathfrak{B}(H) := \left\lbrace (x_k) \in (\mathbb{R}^n)^\mathbb{N}\mid (x_k) \textrm{ is a trajectory of } \eqref{eq:diffinc}\right\rbrace . \]
The set of \textit{feasible} states of the difference inclusion \eqref{eq:diffinc} is the set of states from which a complete trajectory emanates:	
\[ \mathcal{F}(H) := \left\lbrace\xi \mid \exists (x_k)\in \mathfrak{B}(H) \textrm{ with } x_0=\xi \right\rbrace.\] 
The set of \textit{stabilizable} states is the set of states from which a stable trajectory exists:
\[ \mathcal{S}(H) := \{\xi \mid \exists (x_k)\in \mathfrak{B}(H) \textrm{ with } x_0=\xi, \lim_{k\rightarrow\infty} x_k= 0  \}.  \] 
In addition, we define the set of \textit{strongly stable} states as the set of states from which all trajectories are stable:
\[ \bar{\mathcal{S}}(H) := \{\xi \mid \forall (x_k)\in \mathfrak{B}(H) \textrm{ with } x_0=\xi, \lim_{k\rightarrow\infty} x_k= 0  \}.  \] 
It is straightforward to show that if $H$ is a convex process, the sets $\calF(H)$ and $\calS(H)$ are convex cones. The set $\bar{\mathcal{S}}(H)$ is a convex cone if it is not empty. 

We say the system \eqref{eq:diffinc} (or the convex process $H$) is \textit{stabilizable} if every feasible state is stabilizable, that is, $\calF(H)\subseteq\calS(H)$. Similarly, we say it is \textit{strongly stable} if $\calF(H)\subseteq\bar{\mathcal{S}}(H)$. 

The notion of stabilizability considered in this paper differs slightly from the ones employed in e.g. \cite{Phat:96} and \cite{Gajardo:06}, where it is named weak asymptotic stabilizability, and weak asymptotic stability respectively. These papers require all points in $\mathbb{R}^n$ to be stabilizable, which immediately \textit{requires} the process to be strict in order to be stabilizable. This motivates these papers in their assumption of strictness. On the other hand, we argue that in the case of, for example, constrained linear systems, the question whether all feasible states are stabilizable is more natural.

We denote the image of a set $\calS$ under $H$ by $H(\calS) = \{ y\in\mathbb{R}^n \mid \exists x\in \calS \textrm{ s.t. } y\in H(x)\}$. This shorthand allows us to define powers of convex processes, by taking $H^1=H$ and for $q\geq 1$: 
\[H^{q+1}(x) := H(H^{q}(x)) \quad \forall x\in \mathbb{R}^n. \]
We can define the inverse of a convex process by $H^{-1}(y)=\{ x \mid y\in H(x)\}$. Note that this is always defined as a set-valued map. For negative powers of $H$ we use the shorthand $H^{-n}(x) = (H^{-1})^n(x)$.
We define the set of $q$\textit{-step trajectories} as 
\[ \mathfrak{B}_q(H) = \left\lbrace (x_k)_{k=0}^q \in (\mathbb{R}^n)^{q+1}\mid (x_k) \textrm{ satisfies } (\ref{eq:diffinc})\right\rbrace . \] 
Using this, we say that a point $\xi\in\mathbb{R}^n$ is \textit{reachable} if there exists a $q$-step trajectory from the origin, reaching $\xi$. The set of such points is the reachable set:
\[\label{eq:reach def}
\!\!\calR(H)\!=\! \big\lbrace \xi \mid \exists (x_k)_{k=0}^q \in\mathfrak{B}_q(H) \textrm{ s.t. } x_0=0, x_q=\xi \big\rbrace.\]

For a convex cone $ \calC\subseteq\R^n$, we define $\lin(\calC) =-\calC\cap \calC$ and $\Lin(\calC)=\calC-\calC$. It is clear that $\lin(\calC)$ is the largest subspace contained in $\calC$ whereas $\Lin(\calC)$ is the smallest subspace that contains $\calC$. 

Let $H$ be a convex process. Associated with $H$, we define the two linear processes $L_-$ and $L_+$ by
\begin{equation}\label{eq:def of L- and L+}
\graph(L_-)=\lin\big(\gr(H)\big)\textrm{ and } \graph(L_+)=\Lin\big(\graph(H)\big).
\end{equation}
By definition, we therefore have
\begin{equation}\label{eq:lh-+}
\graph(L_-)\subseteq\graph(H)\subseteq\graph(L_+).
\end{equation}
It is clear that $L_-$ and $L_+$ are, respectively, the largest and the smallest (with respect to the graph inclusion) linear processes satisfying \eqref{eq:lh-+}. We call $L_-$ and $L_+$, respectively, the minimal and maximal linear processes associated with $H$. If $H$ is not clear from context, we write $L_-(H)$ and $L_+(H)$ in order to avoid confusion. Similarly, for the respective reachable sets, we write: 
\[ \calR_- = \calR(L_-), \quad \calR_+= \calR(L_+). \] 

For a nonempty set $\calC\subseteq \R^n $, we define the \textit{negative} and \textit{positive polar cone}, respectively,
\begin{align*}
\calC^- &= \{ y\in\R^n \mid \langle x,y\rangle\leq 0 \hspace{1em} \forall x\in \calC\}, \\
\calC^+ &= \{ y\in\R^n \mid \langle x,y\rangle\geq 0 \hspace{1em} \forall x\in \calC\}.
\end{align*}
For a subspace $\calV$, we have that $\calV^-=\calV^+=\calV^\bot$, where the last denotes the orthogonal complement of $\calV$. As such, the positive and negative polar cones generalize the notion of orthogonal complement. Given sets $\calC$ and $\calS$, we have the following:
\begin{equation}\label{eq:dual sum} (\calC+\calS)^-  =\calC^-\cap \calS^-, \quad (\calC\cap\calS)^- = \cl (\calC^-+\calS^-). \end{equation}
The same equations hold for the positive polar cone.

Based on the definition of the negative and positive polar cones, we define \textit{negative} and \textit{positive dual} processes $H^-$ and $H^+$ of $H$ as follows:
\begin{subequations}
\begin{align} \label{eq:def of dual} 
p\in H^-(q) &\iff \langle p,x\rangle\geq \langle q,y\rangle\quad \forall\, (x,y)\in\gr(H),  \\
p\in H^+(q) & \iff \langle p,x\rangle\leq \langle q,y\rangle\quad \forall\, (x,y)\in\gr(H). 
\end{align}
\end{subequations}
The positive dual is alternatively called the \textit{adjoint}, \textit{(left-) transpose} in the literature. 

Note that $H^+(q)= -H^-(-q)$ for all $q$. If $H$ is closed, we know that $(H^+)^- = H$ and
\begin{equation}\label{eq:H(0) to domain} H(0)=(\dom H^+)^+=(\dom H^-)^-.\end{equation} 

If $L$ is a linear process it is clear that the negative and positive dual processes coincide, in which case we denote it by $L^\bot:=L^- =L^+$. The reachable and feasible sets of a linear process $L$ can be determined in a finite number of steps. To be precise $\calF(L)= L^{-n}(\mathbb{R}^n)$ and $\calR(L)= L^n(0)$. The feasible and reachable set of a linear process and its dual are related by:
\vspace{-0.5em}\begin{subequations}\label{eq:LorthRF}
	\begin{gather} 
	\calF(L^\bot) = \calR(L)^\bot, \\
	\calR(L^\bot) = \calF(L)^\bot.
	\end{gather}
\end{subequations}

In addition, the minimal and maximal linear processes associated with a convex process enjoy the following additional properties:
\vspace{-0.5em}\begin{subequations}\label{eq:L(H) vs L(H+-)}
\begin{align}
&L_-(H^-)=L_-(H^+)=L_+^\bot, \\
&L_+(H^-)=L_+(H^+)=L_-^\bot.
\end{align}
\end{subequations}

\begin{lemma}\label{lemm:H and L powers}
	Let $H$ be a convex process and let $L$ be a linear process such that $\graph (L) \subseteq \graph (H)$. For all $x\in\dom(H)$, $y\in\dom(L)$, we have $H(x+y)=H(x)+L(y).$
\end{lemma} 
\BP
Let $x\in\dom(H)$ and $y\in\dom(L)$. We will prove the equality by mutual inclusion. Note that, as $\graph (L) \subseteq \graph (H)$, we know that $L(y)\subseteq H(y)$, and therefore
$$
H(x)+L(y)
\subseteq H(x)+H(y)\subseteq H(x+y).
$$
For the reverse inclusion, first observe that $y\in\dom(L)$ implies that $-y\in\dom(L)$ as $L$ is linear. Then, we have
$$
H(x+y)+L(-y)\subseteq H(x+y)+H(-y)\subseteq H(x).
$$
This shows that $H(x+y)\subseteq H(x)-L(-y)=H(x)+L(y)$, where the last equality follows from $L$ being linear. \EP

The following lemma will show that under a simple condition on $H$, we have that $\calF(H)^- \cap \cl(\calF(H^+))=\{0\}$ and $\calF(H)^- \cap \cl(\calF(H^-))=\{0\}$. These latter conditions will play a central role in the results of Section~\ref{sec:results}.
\begin{lemma}\label{lemm:consequence of dom cond}
	Let $H$ be a convex process such that $\dom H +\calR_-=\mathbb{R}^n$. Then, $\calF(H)^- \cap \cl(\calF(H^+))=\{0\}$ and $\calF(H)^- \cap \cl(\calF(H^-))=\{0\}$. 
\end{lemma}
\BP Let $\dom H +\calR_-=\mathbb{R}^n$. We start with proving the following using induction:
\begin{equation}\label{eq:domHk} \dom H^k +\calR_- =\mathbb{R}^n\quad \forall k\geq 1.\end{equation}
First note that the case of $k=1$ is evident. Now assume that $k\geq 1$ and $\dom H^k +\calR_- =\mathbb{R}^n$. Let $x\in \dom H^k$ and $y\in H^k(x)$. By assumption, we can write $y=\xi-\zeta$, where $\xi\in\dom H$ and $\zeta\in\calR_-$. Since $\calR_-=\calR(L_-)$, we have $\zeta\in L_-^k(\alpha)$ for some $\alpha\in\calR_-$, and therefore we can conclude via Lemma~\ref{lemm:H and L powers} that: 
\[ \xi = y+\zeta\in H^k(x)+L_-^k(\alpha)=H^k(x+\alpha).\]
As $\xi\in\dom H$, this means that $x+\alpha\in \dom H^{k+1}$, and therefore that $\dom H^k \subseteq \dom H^{k+1}+\calR_-$. We can now conclude that
\[ \mathbb{R}^n = \dom H^k+\calR_- \subseteq \dom H^{k+1}+\calR_-,\] 
which concludes the induction step, proving \eqref{eq:domHk}.  

Now let $x\in \dom H^n$ and $y\in H^n(x)$. Again, we can write $y=\xi-\zeta$, where $\xi\in\dom H$ and $\zeta\in\calR_-$. As $\calR_-=\calR(L_-)$, we have that $\zeta\in L_-^n(0)$, and therefore we can use Lemma~\ref{lemm:H and L powers} to conclude that:
\[\xi=y+\zeta \in H^n(x)+L^n_-(0)= H^n(x).\]
Therefore $\dom H^n=\dom H^{n+1}$. By \textit{weak invariance} (see e.g. \cite[Section V]{ReachNullc:19}) this implies that $\dom H^n=\calF(H)$. 

Combining this with \eqref{eq:domHk}, we get that:
\begin{equation}\label{eq:FH+R-=Rn} \mathbb{R}^n=\dom H+\calR_- = \calF(H) +\calR_-. \end{equation} 
By \eqref{eq:L(H) vs L(H+-)} and \eqref{eq:LorthRF}, we have:
\begin{equation}\label{eq:FH+ in R-bot} \calF(H^+) \subseteq \calF(L_+(H^+))= \calF(L_-^\bot)=\calR_-^\bot. \end{equation}
As the right-hand side is closed by definition, we have that $\cl(\calF(H^+))\subseteq \calR_-^\bot$.

Taking the negative polar of \eqref{eq:FH+R-=Rn}, we have by \eqref{eq:dual sum} that $\calF(H)^- \cap \calR_-^\bot=\{0\}$. By the previous, we therefore have that $\calF(H)^- \cap \cl(\calF(H^+))=\zset$. We can repeat the same arguments for $\calF(H)^- \cap \cl\calF(H^-)=\{0\}$. 
\EP

\section{Extended real-valued functions}\label{sec:ervf}

Let $f$ be an extended real-valued function, i.e. a function $f:\mathbb{R}^n\rightarrow \mathbb{R}\cup \{ \pm \infty\}$. We define the \textit{epigraph} of $f$ by
\[ \epi f := \{ (x,\alpha)\in \mathbb{R}^n\times \mathbb{R} \mid \alpha \geq f(x)\}.\]
A function is said to be convex or closed if its epigraph is convex or closed as a set. The \textit{(effective) domain} of the function is the set $\dom f = \{x\in\mathbb{R} \mid f(x)<\infty\}$. The function $f$ is \textit{proper} if $f(x)>-\infty$ for all $x\in\mathbb{R}^n$ and $\dom f\neq \emptyset$.

Let $\calC\subseteq\mathbb{R}^n$. The \textit{indicator function} of $\calC$ is the function $\delta(\cdot \mid \calC)$ given by: 
\[ \delta(x\mid \calC) := \begin{cases} 0 & x\in\calC \\ \infty & x\not\in \calC.\end{cases} \] 
It is straightforward to check that $\dom \delta(\cdot \mid \calC) =\calC$ and that  the indicator function is convex or closed if and only if $\calC$ is.

The function $f$ is called \textit{positively homogeneous of degree} $2$ if $f(\lambda x) = \lambda^2 f(x)$ for all $\lambda \geq 0$ and $x\in \mathbb{R}^n$. If $f$ is positively homogeneous degree 2, it is straightforward to check that the domain of $f$ is a cone. 

A function $f$ is \textit{positive semi-definite} if $f(0)=0$ and $f(x)\geq 0$ for all $x\in \mathbb{R}^n$. Clearly, any such function is proper. Let $\calV$ denote the set of all extended real-valued functions that are closed, convex, positive semi-definite and positively homogeneous of degree $2$. 

%
%

Given an extended real-valued function $f$, we define its \textit{convex conjugate} by 
\[ f^\star (y) := \sup_{x\in\mathbb{R}^n} \left\lbrace y\cdot x -f(x) \right\rbrace. \]
It is well known (see e.g. \cite[Theorem 12.2]{Rockafellar:70}) that the convex conjugate of a proper closed convex function is again a proper closed convex function. In addition, if $f$ is closed and convex, then $f^{\star\star}= f$. As taking the convex conjugate also preserves positive semi-definiteness and positive homogeneity of degree 2, we see that $f\in \calV$ if and only if $f^\star\in \calV$. 

It is straightforward to show that if $f_1(x)\geq f_2(x)$ for all $x\in\mathbb{R}^n$, then $f_1^\star(y)\leq f_2^\star(y)$ for all $y\in \mathbb{R}^n$. 

Given proper convex functions $f$ and $g$, we define their \textit{infimal convolution} $f\ico g$ by the relation $\epi (f\ico g) := \epi f+\epi g$, or equivalently:
\[ (f \ico g) (x) := \inf_{y\in\mathbb{R}^n} f(x-y) +g(y). \] 
A consequence of \cite[Theorem 16.4]{Rockafellar:70} is the following:

\begin{cor}\label{cor: ico conj}
	Let $f$ and $g$ be proper convex functions, then $(f\ico g)^\star= f^\star+g^\star$ and $(\cl (f)+\cl (g))^\star= \cl(f^\star\ico g^\star)$. If in addition $\rint (\dom f)\cap \rint (\dom g) \neq \emptyset$, then $(f+g)^\star= f^\star \ico g^\star$. 
\end{cor} 

\begin{example}\label{ex:conv conj}
	Let $\calC$ be a convex cone, then the following are straightforward to check:
	\vspace{-0.5em}\begin{align*} 
	f(x) = \tfrac{1}{2}\norm{x}^2 &\implies f^\star(y) = \tfrac{1}{2}\norm{y}^2, \\
	g(x)=\delta(x\mid \calC) &\implies g^\star(y) = \delta(y\mid \calC^-).
	\end{align*}	
	Denote the squared \textit{distance function} to $\calC$ as
	\vspace{-0.5em}\[ h(x) := \inf_{y\in \calC} \tfrac{1}{2}\norm{x-y}^2 = (f \ico g) (x).\] 
	We can use Corollary~\ref{cor: ico conj} to state: 
	\vspace{-0.5em}\[ h^\star(y) = \tfrac{1}{2}\norm{y}^2+\delta(y\mid \calC^-). \]
\end{example}

We define the \textit{restriction of} $f$ \textit{to} $\calC$ by
\[ f|_\calC(x)  = f(x)+\delta(x\mid \calC). \]
If $f\in\calV$ then for any closed convex cone $\calC$ we have that $f|_\calC\in\calV$. 

\begin{definition}\label{def:posdef}
	Let $\calC$ be a convex cone. We say a function $f\in\calV$ is \textit{positive definite with respect to} $\calC$ if there exist $0<\alpha\leq\beta <\infty$ such that for all $x\in \calC$:
	\[ \alpha\norm{x}^2\leq f(x) \leq \beta\norm{x}^2.\]

\end{definition}
In particular, this implies that if $f$ is positive definite with respect to $\calC$, then $0<f(x)<\infty$ for all $x\in \calC$. Clearly, this also implies that $\calC\subseteq \dom f $. 

It is well known that if $f$ is positive definite with respect to $\mathbb{R}^n$, then so is $f^\star$. This result is generalized in the following theorem. 

\begin{theorem}\label{lemm: conds for conj in V} 
	Let $f\in\calV$ and let $\calC,\calD\subseteq\mathbb{R}^n$ be convex cones such that $f$ is positive definite with respect to $\calC$ and $\calC^-\cap \cl(\calD)= \{0\}$. Then, $(f|_\calC)^\star\in\calV$ and $(f|_\calC)^\star$ is positive definite with respect to $\calD$. 
\end{theorem} 
\BP First, we will prove that $(f|_{\calC})^\star=(f|_{\cl(\calC)})^\star$. Note that $\rint(\calC)\cap \rint(\dom f)=\rint(\calC)$, and therefore we know by Corollary~\ref{cor: ico conj} that:
\begin{align*} 
(f|_{\calC})^\star &= (f+\delta(\cdot\mid \calC))^\star \\
&= f^\star \ico \delta(\cdot \mid \calC^-). 
\end{align*}
As the convex conjugate of a function is always closed, and therefore equal to its closure, we can apply Corollary~\ref{cor: ico conj} again, leading to: 
\begin{align*} 
(f|_{\calC})^\star	&= \cl \left(f^\star \ico \delta(\cdot \mid \calC^-\right) \\
&= \left(\cl(f) + \cl(\delta(\cdot\mid \calC))\right)^\star \\
&= (f|_{\cl(\calC)})^\star.
\end{align*}
Then, as $\cl(\calC)$ is closed, we have $f|_{\cl(\calC)}\in \calV$ and therefore we know that $(f|_{\calC})^\star\in\calV$. 

From Definition~\ref{def:posdef}, we know there exists $0<\alpha<\infty$ such that
\[\alpha\norm{x}^2+\delta(x\mid \calC)\leq f|_\calC(x).\]
Recall that $f|_\calC(x)\geq h(x)$ for all $x\in\mathbb{R}^n$ implies that $(f|_\calC)^\star(y)\leq h^\star(y)$ for all $y\in \mathbb{R}^n$. By using Corollary~\ref{cor: ico conj} and Example~\ref{ex:conv conj} we therefore have that for all $y\in\mathbb{R}^n$:
\[ (f|_\calC)^\star(y) \leq \frac{1}{4\alpha} \inf_{x\in \calC^-} \norm{y-x}^2.\] 
Note that $0\in\calC^-$, and therefore:
\begin{equation}\label{eq:one side} (f|_\calC)^\star(y) \leq \frac{1}{4\alpha} \norm{y}^2. \end{equation}

Before we are able to prove that $(f|_\calC)^\star$ is positive definite with respect to $\calD$, we require a preliminary result. Given $\calD$, we define the function $h$ by:
\[h(y) :=  \inf_{x\in \calC^-} \norm{y-x}^2 +\delta(y\mid \calD).\]
Taking $\calS = \{ y \mid \norm{y}^2=1\}$, we have that: 
\begin{align*} \inf_{y\in\calS} h(y)  &=\inf_{y\in \calS\cap\calD}  \inf_{x\in \calC^-} \norm{y-x}^2 \\
				&\geq \inf_{y\in \calS\cap\cl(\calD)}  \inf_{x\in \calC^-} \norm{y-x}^2 .\end{align*}
We will now prove by contradiction that the last term of the previous inequality is strictly greater than $0$. Assume the contrary, i.e. that there exist sequences $(x_k)_{k=0}^\infty,(y_k)_{k=0}^\infty$ where $x_k\in \calC^-$ and $y_k\in \calS\cap\cl(\calD)$ for all $k\geq0$, such that $\lim_{k\rightarrow \infty} \norm{y_k-x_k}^2=0$. 

Note that $\calS\cap\cl(\calD)$ is both closed and bounded. Therefore there exists a subsequence $(y_{k_\ell})_{\ell=0}^\infty$ that converges to $\bar{y}\in \calS\cap\cl(\calD)$. We can now use the triangle inequality to show that
\[\lim_{\ell\rightarrow\infty} \norm{\bar{y}-x_{k_\ell}}^2 \leq \lim_{\ell\rightarrow\infty} \norm{\bar{y}-y_{k_\ell}}^2+ \lim_{\ell\rightarrow\infty}\norm{y_{k_\ell}-x_{k_\ell}}^2 =0. \]
This means that $\bar{y}$ is a limit point of $\calC^-$. As $\calC^-$ is closed, this implies that $\bar{y}\in \calC^-\cap(\calS\cap\cl(\calD))$. Recall that by assumption $\calC^-\cap \cl(\calD)=\zset$ and $0\not\in\calS$, thus leading to a contradiction.  

Following the previous, we know that there exists $\gamma>0$ such that $\inf_{y\in\calS} h(y) \geq \gamma$. As $h$ is positively homogeneous of degree 2, we can conclude that $h(y) \geq \gamma \norm{y}^2$ for all $y\in\mathbb{R}^n$.

Similar to the first part of this proof, there exists $0<\beta<\infty$ such that for each $y\in\mathbb{R}^n$ we have the inequality:
\[ \frac{1}{4\beta} \inf_{x\in \calC^-} \norm{y-x}^2 \leq (f|_\calC)^\star(y).\]
Using our previous results, we therefore have that:
\[ (f|_\calC)^\star(y)+\delta(y\mid \calD) \geq  \frac{1}{4\beta} h(y) \geq \frac{\gamma}{4\beta} \norm{y}^2.\]
Combining this with \eqref{eq:one side} shows that $(f|_\calC)^\star$ is positive definite with respect to $\calD$, thus proving the theorem. \EP

\begin{remark} 
	Indeed, if $f$ is positive definite with respect to $\mathbb{R}^n$, we can take $\calD= \mathbb{R}^n$ and see that $f^\star=(f|_{\mathbb{R}^n})^\star\in \calV$ is positive definite with respect to $\mathbb{R}^n$. 
\end{remark}
\section{Lyapunov functions}\label{sec:Lyap}
These preliminaries lead us to a definition of Lyapunov functions for convex processes. 
\begin{definition}\label{def:Lyap}
	Let $H$ be a convex process. A function $V\in \calV$ is a \textit{weak Lyapunov function for} $H$ if $V$ is positive definite with respect to $\calF(H)$ and there exists $\gamma \in (0,1)$ such that 
	\begin{equation}\label{eq: weak lyap} 
	\forall x\in \calF(H), \enskip \exists y\in \calF(H)\cap  H(x) \textrm{ s.t. }  V(y) \leq \gamma V(x). \end{equation} 
	A function $V\in \calV$ is a \textit{strong Lyapunov function for} $H$ if $V$ is positive definite with respect to $\calF(H)$ and there exists $\gamma \in (0,1)$ such that 
	\begin{equation}\label{eq: strong lyap} 
	\forall x\in \calF(H), \enskip \forall y\in \calF(H)\cap H(x) \quad V(y) \leq \gamma V(x). \end{equation} 	
\end{definition}

We will first compare this definition to earlier notions of Lyapunov functions for convex processes. In terms of the notation of this paper, the definition for a weak Lyapunov function used in \cite{Goebel:13} requires $V\in\calV$ to be positive definite with respect to $\mathbb{R}^n$ and the existence of $\gamma\in(0,1)$ such that
\begin{equation}\label{eq:Goebel wL}	\forall x\in \dom H, \enskip \exists y\in  H(x) \textrm{ s.t. }  V(y) \leq \gamma V(x). \end{equation}
Similarly, (strong) Lyapunov functions are defined in \cite{Goebel:13} as functions $V\in\calV$ that are positive definite with respect to $\mathbb{R}^n$ such that there exists $\gamma\in(0,1)$ for which 
\begin{equation}\label{eq:Goebel sL}	\forall x\in \dom H, \enskip \forall y\in  H(x) \quad V(y) \leq \gamma V(x). \end{equation}

Note that for strict convex processes $\dom H = \calF(H)=\mathbb{R}^n$, which makes these definitions coincide with Definition~\ref{def:Lyap}. On the other hand, for nonstrict convex processes, important differences arise. Using the following two examples, we will argue that Definition~\ref{def:Lyap} is more natural for both weak and strong Lyapunov functions. First, we consider weak Lyapunov functions. The following is an example of a convex process and a function that is a weak Lyapunov function in the sense of \eqref{eq:Goebel wL}, which fails to be stabilizable. 
\begin{example} 
	Let $H:\mathbb{R}^2\rightrightarrows \mathbb{R}^2$ be the convex process given by:
	\[ H(x) = \begin{bmatrix} 0 & 0 \\ 0 & \tfrac{1}{2} \end{bmatrix} x + \mathbb{R}\times \{0\}\quad  \textrm{if } \begin{bmatrix} 1 & -2 \\ 0 & 1 \end{bmatrix}x\geq 0, \] 
	and empty otherwise. Here `$\geq$' is understood to hold element-wise. 
%
	It is straightforward to check that $\calF(H)=\dom H$. Let $V(x) = \tfrac{1}{2}\norm{x}^2$, then for each $x\in \dom H$, there exists $y\in H(x)$ such that $V(y) \leq \tfrac{1}{4} \norm{x}^2$. On the other hand, it is straightforward to check that $H$ is not stabilizable. Indeed, for $x\in\calF(H)$, we have that $y\in\calF(H)\cap H(x)$ implies that $V(y)\geq V(x)$. 
\end{example}

The following is an example of a convex process which is strongly stable, but for which there does not exist any strong Lyapunov function in the sense of \eqref{eq:Goebel sL}. 
\begin{example}
	Let $H:\mathbb{R}^2\rightrightarrows \mathbb{R}^2$ be the convex process given by:
	\[ H(x) = -\tfrac{1}{2}x + \{0\} \times \mathbb{R}_- \quad \textrm{if } \begin{bmatrix} 0 & 1 \end{bmatrix} x \geq 0, \] 
	and empty for other $x$. 
	It is straightforward to see that $\calF(H)= \{ x \mid \begin{bmatrix} 0 & 1 \end{bmatrix} x = 0 \}$, and that $H(x)\cap \calF(H)$ contains only $-\tfrac{1}{2}x$ for any $x\in \calF(H)$. Therefore $H$ is stable. However, for e.g. the point $x=\left(\!\begin{smallmatrix} 0 \\ 1\end{smallmatrix}\!\right)$, we have that $x\in \dom H$ and $H(x) = \{ \left(\!\begin{smallmatrix} 0 \\ \alpha \end{smallmatrix}\!\right) \mid \alpha \leq -\tfrac{1}{2}\}$. Clearly, for this $x$, there does not exists a function $V\in\calV$ that is positive definite with respect to $\mathbb{R}^n$ and $\gamma\in (0,1)$ for which we have that $V(y)\leq \gamma V(x)$ for each $y\in H(x)$. 
\end{example} 

\section{Duality theorems}\label{sec:results}
A direct consequence of Definition~\ref{def:Lyap} arises when considering strong Lyapunov functions. Let $H$ be a convex process. If there exists a strong Lyapunov function $V$ for $H$, then for every $y\in\calF(H) \cap H(0)$ we have $V(y)\leq \gamma V(0)$. As $y\in\calF(H)$ and $V$ is positive definite with respect to $\calF(H)$, this implies that $y=0$. This means that a necessary condition for the existence of a strong Lyapunov function is that $\calF(H) \cap H(0)=\zset$.

The following theorem will relate weak Lyapunov functions for $H$ to strong Lyapunov functions for $H^+$. By the previous and the equality from \eqref{eq:H(0) to domain}, $H^+$ admits a strong Lyapunov function only if $\calF(H^+)\cap(\dom H)^- =\zset$. 

In fact, the following theorem works under a slightly stronger assumption, namely that $\calF(H)^- \cap \cl(\calF(H^+))=\{0\}$. Note that by Lemma~\ref{lemm:consequence of dom cond} we can guarantee that this last condition holds if $\dom H +\calR_-=\mathbb{R}^n$.

\begin{theorem}\label{thm: weak-strong}
	 Let $H$ be a closed convex process such that $\calF(H)^- \cap \cl(\calF(H^+))=\{0\}$. If $V\in \calV$ is a weak Lyapunov function for $H$, then $W:=(V|_{\calF(H)})^\star$ is a strong Lyapunov function for $H^+$. 
\end{theorem}

\BP We can use Theorem~\ref{lemm: conds for conj in V} to conclude that $W\in\calV$ and $W$ is positive definite with respect to $\calF(H^+)$. What remains is to prove that there exists $\gamma \in (0,1)$ such that \eqref{eq: strong lyap} holds for $H^+$ and $W$. Let $q\in \calF(H^+)$ and $p\in \calF(H^+)\cap H^+(y)$, then 
\begin{align*}
W(p) &= \left(V+\delta(\cdot \mid \calF(H))\right)^\star(p) \\
&= \sup_{x\in \calF(H)} \{p\cdot x -V(x)\}. \end{align*}

For $x\in \calF(H)$, we know that $V(x)\geq \tfrac{1}{\gamma}V(y)$ for some $y\in \calF(H)\cap H(x)$, therefore 
\begin{align*} 
W(p)	& \leq  \sup_{x\in \calF(H)} \{p\cdot x -\frac{1}{\gamma}\inf_{y\in \calF(H)\cap H(x)}V(y) \} \\
		&= \sup_{x\in \calF(H)}\sup_{y\in \calF(H)\cap H(x)} \{ p\cdot x -\frac{1}{\gamma}V(y) \}. 
\end{align*}
Following the definition of $H^+$, we see that $p\cdot x\leq q\cdot y$:
\[ W(p) \leq \sup_{x\in \calF(H)}\sup_{y\in \calF(H)\cap H(x)} \{ q\cdot y -\frac{1}{\gamma}V(y)\}  \]
Note that for each $x\in\calF(H)$ the set $\calF(H)\cap H(x)$ is nonempty by definition, and therefore
\begin{align*} 
W(p)	&\leq \frac{1}{\gamma}  \sup_{y\in \calF(H)} \{ \gamma q\cdot y -V(y)\} \\
		&= \frac{1}{\gamma} W(\gamma q) = \gamma W(q). 
\end{align*}
Thus proving that \eqref{eq: strong lyap} holds, which proves the theorem.  
\EP 

Note that if $H$ is a strict convex process, the condition $\dom H+\calR_-=\mathbb{R}^n$ holds immediately. This means that the previous theorem is a generalization of \cite[Theorem 2.2]{Goebel:13}. 

The following theorem, which is based on \cite[Theorem  2.4]{Goebel:13}, provides conditions under which a strong Lyapunov function for a convex process $H$ can be transformed into a weak Lyapunov function for another convex process. 

\begin{theorem} \label{thm: strong-weak}
	 Let $H$ and $G$ be closed convex processes such that $\calF(H)^- \cap \cl (\calF(G))=\{0\}$ and for all $ x\in \calF(H)$ and $q\in \calF(G)$:
	 \begin{equation}\label{eq:assumption} \inf_{ p\in \calF(G)\cap G(q)} p\cdot x \leq \sup_{y\in \calF(H) \cap H(x)} y\cdot q. \end{equation}	 
	  If $V\in \calV$ is a strong Lyapunov function for $H$, then $W:=(V|_{\calF(H)})^\star$ is a weak Lyapunov function for $G$. 
\end{theorem}
\BP We can follow the proof of Theorem~\ref{thm: weak-strong}, to see that $W\in\calV$ and $W$ is positive definite with respect to $\calF(G)$. 

Note that to prove the theorem it suffices to prove that for every $q\in \calF(G)$, the following holds: 
\[ \inf_{p\in \calF(G)\cap G(q)} W(p) \leq \gamma W(q). \] 
Let $q\in \calF(G)$, then:
\[\inf_{p\in \calF(G)\cap G(q)} W(p) = \inf_{p\in \calF(G)\cap G(q)} \sup_{x\in\calF(H)} \{ p\cdot x -V(x)\}.\]

As $V|_{\calF(H)}$ is coercive, we can swap the infimum and supremum using \cite[Theorem 37.3]{Rockafellar:70}, leading to:
\[\inf_{p\in \calF(G)\cap G(q)} W(p) =  \sup_{x\in\calF(H)}\inf_{p\in \calF(G)\cap G(q)} \{ p\cdot x -V(x)\}.\]
Then, we can apply the inequality of \eqref{eq:assumption} to obtain: 
\[\inf_{p\in \calF(G)\cap G(q)} W(p) \leq \sup_{x\in\calF(H)}\sup_{y\in \calF(H)\cap H(x)}  \{ y\cdot q -V(x)\}.\]
As $V$ is a strong Lyapunov function for $H$, we know that $V(x)\geq \tfrac{1}{\gamma}V(y)$ for all $y\in \calF(H)\cap H(x)$. 
\begin{align*} \inf_{p\in \calF(G)\cap G(q)} W(p) &\leq \sup_{x\in\calF(H)}\sup_{y\in \calF(H)\cap H(x)}  \{ y\cdot q -\tfrac{1}{\gamma}V(y)\}\\
&\leq \sup_{y\in \calF(H)}  \{ y\cdot q -\tfrac{1}{\gamma}V(y)\} = \gamma W(q).\end{align*}
Thus proving the theorem. \EP

So far, checking the conditions under which the previous theorem works might seem daunting. However, for a few choices of $G$, we can check this condition easily. 
\begin{example} 
	If $G=H^+$, we see that the assumption that \eqref{eq:assumption} is satisfied for all $ x\in \calF(H)$ and $q\in \calF(H^+)$ holds immediately. As any strong Lyapunov function is a weak Lyapunov function, this reduces Theorem~\ref{thm: strong-weak} to a specific case of Theorem~\ref{thm: weak-strong}.
\end{example} 
\begin{example}
	If $G=H^-$, the assumption that \eqref{eq:assumption} is satisfied for all $x\in \calF(H)$ and $q\in \calF(H^-)$ is equivalent to the assumption that for all $x\in \calF(H)$ and $q\in \calF(H^+)$:
	\[ \sup_{ p\in \calF(H^+)\cap H^+(q)} p\cdot x \leq \inf_{y\in \calF(H) \cap H(x)} y\cdot q. \]
	In turn, we can use \cite[Theorem 2.9]{Smirnov:02} to show that this follows if for example $\calF(H) \subseteq \rint \dom H$. As an example, for linear processes this condition is always satisfied.
\end{example}
\section{Conclusions}\label{sec:conc}
\vspace{-0.3em}In this paper, we have provided a new definition for Lyapunov functions for difference inclusions of nonstrict convex processes. As shown in a few examples, this definition better captures the stabilizability properties of these systems. Building on this definition, we have shown that under certain conditions, a weak Lyapunov functions for a convex process can naturally be transformed to a strong Lyapunov function for its dual. In addition we reveal conditions under which a strong Lyapunov function can be transformed to a weak Lyapunov function for another associated convex process. These results generalize known results and the conditions required on the domain are, in some sense, close to being necessary. 
\vspace{-1em}
\bibliography{convex_processes}
\bibliographystyle{IEEEtran}

\end{document}

%% file: ka-newcommands-siam.tex
\DeclareMathOperator{\im}{im}

\DeclareMathOperator{\cl}{cl}
\DeclareMathOperator{\rint}{ri}

\DeclareMathOperator{\epi}{epi}
\DeclareMathOperator{\dom}{dom}
\DeclareMathOperator{\graph}{gr}
\DeclareMathOperator{\gr}{gr}
\DeclareMathOperator{\lin}{lin}
\DeclareMathOperator{\Lin}{Lin}

\let\leq\leqslant
\let\geq\geqslant
\let\emptyset\varnothing

\newcommand{\calC}{\ensuremath{\mathcal{C}}}
\newcommand{\calD}{\ensuremath{\mathcal{D}}}

\newcommand{\calF}{\ensuremath{\mathcal{F}}}

\newcommand{\calR}{\ensuremath{\mathcal{R}}}
\newcommand{\calS}{\ensuremath{\mathcal{S}}}
\newcommand{\calT}{\ensuremath{\mathcal{T}}}

\newcommand{\calV}{\ensuremath{\mathcal{V}}}

\newcommand{\bmat}{\begin{matrix}}
\newcommand{\emat}{\end{matrix}}
\newcommand{\bbm}{\begin{bmatrix}}
\newcommand{\ebm}{\end{bmatrix}}
\newcommand{\bbma}{\begin{bmatrix*}[r]}
\newcommand{\ebma}{\end{bmatrix*}}
\newcommand{\bpm}{\begin{pmatrix}}
\newcommand{\epm}{\end{pmatrix}}
\newcommand{\bpma}{\begin{pmatrix*}[r]}
\newcommand{\epma}{\end{pmatrix*}}
\newcommand{\bvm}{\begin{vmatrix}}
\newcommand{\evm}{\end{vmatrix}}
\newcommand{\bse}{\begin{subequations}}
\newcommand{\ese}{\end{subequations}}
\newcommand{\beq}{\begin{equation}}
\newcommand{\eeq}{\end{equation}}

\newcommand{\ben}{\begin{enumerate}}
\newcommand{\een}{\end{enumerate}}
\newcommand{\beni}{\begin{enumerate}}
\newcommand{\eeni}{\end{enumerate}}
\newcommand{\bena}{\begin{enumerate}}
\newcommand{\eena}{\end{enumerate}}

\newcommand{\bit}{\begin{itemize}}
\newcommand{\eit}{\end{itemize}}
\newcommand{\bthe}{\begin{theorem}}
\newcommand{\ethe}{\end{theorem}}
\newcommand{\blem}{\begin{lemma}}
\newcommand{\elem}{\end{lemma}}
\newcommand{\bprop}{\begin{proposition}}
\newcommand{\eprop}{\end{proposition}}
\newcommand{\bex}{\begin{example}}
\newcommand{\eex}{\end{example}}
\newcommand{\bas}{\begin{assumption}}
\newcommand{\eas}{\end{assumption}}
\newcommand{\bre}{\begin{remark}}
\newcommand{\ere}{\end{remark}}
\newcommand{\bcor}{\begin{corollary}}
\newcommand{\ecor}{\end{corollary}}
\newcommand{\bdfn}{\begin{definition}}
\newcommand{\edfn}{\end{definition}}
\newcommand{\bcon}{\begin{conjecture}}
\newcommand{\econ}{\end{conjecture}}

\newcommand{\pset}[1]{\ensuremath{\{#1\}}}

\newcommand{\zset}{\ensuremath{\pset{0}}}

\newcommand{\norm}[1]{\ensuremath{\| #1 \|}}

\newcommand{\R}{\ensuremath{\mathbb R}}

\newcommand{\N}{\ensuremath{\mathbb N}}

\newcommand{\BP}{\noindent{\bf Proof. }}
\newcommand{\EP}{\hspace*{\fill} $\blacksquare$\noindent}

%% file: 2020-Lyap.bbl
\begin{thebibliography}{10}
\providecommand{\url}[1]{#1}
\csname url@samestyle\endcsname
\providecommand{\newblock}{\relax}
\providecommand{\bibinfo}[2]{#2}
\providecommand{\BIBentrySTDinterwordspacing}{\spaceskip=0pt\relax}
\providecommand{\BIBentryALTinterwordstretchfactor}{4}
\providecommand{\BIBentryALTinterwordspacing}{\spaceskip=\fontdimen2\font plus
\BIBentryALTinterwordstretchfactor\fontdimen3\font minus
  \fontdimen4\font\relax}
\providecommand{\BIBforeignlanguage}[2]{{%
\expandafter\ifx\csname l@#1\endcsname\relax
\typeout{** WARNING: IEEEtran.bst: No hyphenation pattern has been}%
\typeout{** loaded for the language `#1'. Using the pattern for}%
\typeout{** the default language instead.}%
\else
\language=\csname l@#1\endcsname
\fi
#2}}
\providecommand{\BIBdecl}{\relax}
\BIBdecl

\bibitem{Rockafellar:67}
R.~T. Rockafellar, \emph{Monotone Processes of Convex and Concave Type}, ser.
  American Mathematical Society memoirs.\hskip 1em plus 0.5em minus 0.4em\relax
  American Mathematical Society, 1967.

\bibitem{Rockafellar:70}
------, \emph{Convex Analysis}, ser. Princeton Mathematical Series, No.
  28.\hskip 1em plus 0.5em minus 0.4em\relax Princeton, N.J.: Princeton
  University Press, 1970.

\bibitem{makarov:77}
V.~L. Makarov and A.~M. Rubinov, \emph{Mathematical Theory of Economic Dynamics
  and Equilibria}.\hskip 1em plus 0.5em minus 0.4em\relax Springer-Verlag,
  1977.

\bibitem{j37}
M.~D. Kaba and M.~K. Camlibel, ``A spectral characterization of controllability
  for linear discrete-time systems with conic constraints,'' \emph{SIAM Journal
  on Control and Optimization}, vol.~53, no.~4, pp. 2350--2372, 2015.

\bibitem{oh:05-2}
S.~R. Oh and S.~K. Agrawal, ``A reference governor-based controller for a cable
  robot under input constraints,'' \emph{IEEE Transactions on Control Systems
  Technology}, vol.~13, no.~4, pp. 639--645, 2005.

\bibitem{Angeli:09}
D.~Angeli, P.~De~Leenheer, and E.~D. Sontag, ``Chemical networks with inflows
  and outflows: A positive linear differential inclusions approach,''
  \emph{Biotechnology Progress}, vol.~25, no.~3, pp. 632--642, 2009.

\bibitem{Frankowska:87}
H.~Frankowska, ``Local controllability and infinitesimal generators of
  semigroups of set-valued maps,'' \emph{SIAM Journal on Control and
  Optimization}, vol.~25, no.~2, pp. 412--432, 1987.

\bibitem{Aubin:84}
J.-P. Aubin and I.~Ekeland, \emph{Applied nonlinear analysis}, ser. Pure and
  applied mathematics.\hskip 1em plus 0.5em minus 0.4em\relax J. Wiley, 1984.

\bibitem{AFO:86}
J.-P. Aubin, H.~Frankowska, and C.~Olech, ``Controllability of convex
  processes,'' \emph{SIAM Journal on Control and Optimization}, vol.~24, no.~6,
  pp. 1192--1211, 1986.

\bibitem{PD:94}
V.~N. Phat and T.~C. Dieu, ``On the {K}re\u\i n-{R}utman theorem and its
  applications to controllability,'' \emph{Proceedings of the American
  Mathematical Society}, vol. 120, no.~2, pp. 495--500, 1994.

\bibitem{Smirnov:02}
G.~V. Smirnov, \emph{Introduction to the Theory of Differential Inclusions},
  ser. Graduate Studies in Mathematics.\hskip 1em plus 0.5em minus 0.4em\relax
  Rhode Island: American Mathematical Society, 2002, vol.~41.

\bibitem{Phat:96}
V.~N. Phat, ``Weak asymptotic stabilizability of discrete-time systems given by
  set-valued operators,'' \emph{Journal of Mathematical Analysis and
  Applications}, vol. 202, no.~2, pp. 363 -- 378, 1996.

\bibitem{Seeger:01}
A.~Seeger, ``A duality result for the reachable cone of a nonstrict convex
  process,'' \emph{Journal of Nonlinear and Convex Analysis}, vol.~2, no.~3,
  pp. 363--368, 2001.

\bibitem{ReachNullc:19}
J.~Eising and M.~K. Camlibel, ``On reachability and null-controllability of
  nonstrict convex processes,'' \emph{IEEE Control Systems Letters}, vol.~3,
  no.~3, pp. 751--756, 2019.

\bibitem{Kellet:04a}
C.~M. Kellett and A.~R. Teel, ``Discrete-time asymptotic controllability
  implies smooth control-{Lyapunov} function,'' \emph{Systems \& Control
  Letters"}, vol.~52, no.~5, pp. 349 -- 359, 2004.

\bibitem{Kellet:04}
------, ``Smooth {Lyapunov} functions and robustness of stability for
  difference inclusions,'' \emph{Systems \& Control Letters}, vol.~52, no.~5,
  pp. 395 -- 405, 2004.

\bibitem{Goebel:11}
R.~Goebel, ``Set-valued {Lyapunov} functions for difference inclusions,''
  \emph{Automatica}, vol.~47, no.~1, pp. 127 -- 132, 2011.

\bibitem{Goebel:13}
------, ``Lyapunov functions and duality for convex processes,'' \emph{SIAM
  Journal on Control and Optimization}, vol.~51, pp. 3332--3350, 2013.

\bibitem{Willems:91}
J.~C. Willems, ``Paradigms and puzzles in the theory of dynamical systems,''
  \emph{IEEE Transactions on Automatic Control}, vol.~36, no.~3, pp. 259--294,
  1991.

\bibitem{Gajardo:06}
P.~Gajardo and A.~Seeger, ``Higher-order spectral analysis and weak asymptotic
  stability of convex processes,'' \emph{Journal of Mathematical Analysis and
  Applications}, vol. 318, no.~1, pp. 155--174, 2006.

\end{thebibliography}
